\documentstyle[11pt]{article}
\begin{document}
%
%
%
\newtheorem{theorem}      {Th\'eor\`eme}[section]
\newtheorem{theorem*}     {theorem}
\newtheorem{proposition}  [theorem]{Proposition}
\newtheorem{definition}   [theorem]{Definition}
\newtheorem{e-lemme}        [theorem]{Lemma}
\newtheorem{cor}   [theorem]{Corollaire}
\newtheorem{resultat}     [theorem]{R\'esultat}
\newtheorem{eexercice}    [theorem]{Exercice}
\newtheorem{rrem}    [theorem]{Remarque}
\newtheorem{pprobleme}    [theorem]{Probl\`eme}
\newtheorem{eexemple}     [theorem]{Exemple}
\newcommand{\preuve}      {\paragraph{Preuve}}
\newenvironment{probleme} {\begin{pprobleme}\rm}{\end{pprobleme}}
\newenvironment{remarque} {\begin{rremarque}\rm}{\end{rremarque}}
\newenvironment{exercice} {\begin{eexercice}\rm}{\end{eexercice}}
\newenvironment{exemple}  {\begin{eexemple}\rm}{\end{eexemple}}
%
%
\newtheorem{e-theo}      [theorem]{Theorem}
\newtheorem{theo*}     [theorem]{Theorem}
\newtheorem{e-pro}  [theorem]{Proposition}
\newtheorem{e-def}   [theorem]{Definition}
\newtheorem{e-lem}        [theorem]{Lemma}
\newtheorem{e-cor}   [theorem]{Corollary}
\newtheorem{e-resultat}     [theorem]{Result}
\newtheorem{ex}    [theorem]{Exercise}
\newtheorem{e-rem}    [theorem]{Remark}
\newtheorem{prob}    [theorem]{Problem}
\newtheorem{example}     [theorem]{Example}
\newcommand{\proof}         {\paragraph{Proof~: }}
\newcommand{\hint}          {\paragraph{Hint}}
\newcommand{\heuristicproof}{\paragraph{heuristic proof}}
\newenvironment{e-probleme} {\begin{e-pprobleme}\rm}{\end{e-pprobleme}}
\newenvironment{e-remarque} {\begin{e-rremarque}\rm}{\end{e-rremarque}}
\newenvironment{e-exercice} {\begin{e-eexercice}\rm}{\end{e-eexercice}}
\newenvironment{e-exemple}  {\begin{e-eexemple}\rm}{\end{e-eexemple}}
\newcommand{\reell}    {{{\rm I\! R}^l}}
\newcommand{\reeln}    {{{\rm I\! R}^n}}
\newcommand{\reelk}    {{{\rm I\! R}^k}}
\newcommand{\reelm}    {{{\rm I\! R}^m}}
\newcommand{\reelp}    {{{\rm I\! R}^p}}
\newcommand{\reeld}    {{{\rm I\! R}^d}}
\newcommand{\reeldd}   {{{\rm I\! R}^{d\times d}}}
\newcommand{\reelnn}   {{{\rm I\! R}^{n\times n}}}
\newcommand{\reelnd}   {{{\rm I\! R}^{n\times d}}}
\newcommand{\reeldn}   {{{\rm I\! R}^{d\times n}}}
\newcommand{\reelkd}   {{{\rm I\! R}^{k\times d}}}
\newcommand{\reelkl}   {{{\rm I\! R}^{k\times l}}}
\newcommand{\reelN}    {{{\rm I\! R}^N}}
\newcommand{\reelM}    {{{\rm I\! R}^M}}
\newcommand{\reelplus} {{{\rm I\! R}^+}}
\newcommand{\reelo}    {{{\rm I\! R}\setminus\{0\}}}
\newcommand{\reld}    {{{\rm I\! R}_d}}
\newcommand{\relplus} {{{\rm I\! R}_+}}
\newcommand{\1}        {{\bf 1}}

\newcommand{\cov}      {{\hbox{cov}}}
\newcommand{\sss}      {{\cal S}}
\newcommand{\indic}    {{{\rm I\!\! I}}}
\newcommand{\pp}       {{{\rm I\!\!\! P}}}
\newcommand{\qq}       {{{\rm I\!\!\! Q}}}
\newcommand{\ee}       {{{\rm I\! E}}}

\newcommand{\B}        {{{\rm I\! B}}}
\newcommand{\cc}       {{{\rm I\!\!\! C}}}
\renewcommand{\pp}        {{{\rm I\!\!\! P}}}
\newcommand{\HHH}        {{{\rm I\! H}}}
\newcommand{\N}        {{{\rm I\! N}}}
\newcommand{\R}        {{{\rm I\! R}}}
\newcommand{\D}        {{{\rm I\! D}}}
\newcommand{\Z}       {{{\rm Z\!\! Z}}}
\newcommand{\C}        {{\bf C}}        
\newcommand{\T}        {{\bf T}}        
\newcommand{\E}        {{\bf E}}        
\newcommand{\rfr}[1]    {\stackrel{\circ}{#1}}
\newcommand{\equiva}    {\displaystyle\mathop{\simeq}}
\newcommand{\eqdef}     {\stackrel{\triangle}{=}}
\newcommand{\limps}     {\mathop{\hbox{\rm lim--p.s.}}}
\newcommand{\Limsup}    {\mathop{\overline{\rm lim}}}
\newcommand{\Liminf}    {\mathop{\underline{\rm lim}}}
\newcommand{\Inf}       {\mathop{\rm Inf}}
\newcommand{\vers}      {\mathop{\;{\rightarrow}\;}}
\newcommand{\versup}    {\mathop{\;{\nearrow}\;}}
\newcommand{\versdown}  {\mathop{\;{\searrow}\;}}
\newcommand{\vvers}     {\mathop{\;{\longrightarrow}\;}}
\newcommand{\cvetroite} {\mathop{\;{\Longrightarrow}\;}}
\newcommand{\ieme}      {\hbox{i}^{\hbox{\smalltype\`eme}}}
\newcommand{\eqps}      {\, \buildrel \rm \hbox{\rm\smalltype p.s.} \over =
\,}
\newcommand{\eqas}      {\,\buildrel\rm\hbox{\rm\smalltype a.s.} \over = \,}
\newcommand{\argmax}    {\hbox{{\rm Arg}}\max}
\newcommand{\argmin}    {\hbox{{\rm Arg}}\min}
\newcommand{\indep}{\perp\!\!\!\!\perp}
\newcommand{\abs}[1]{\left| #1 \right|}
\newcommand{\crochet}[2]{\langle #1 \,,\, #2 \rangle}
\newcommand{\espc}[3]   {E_{#1}\left(\left. #2 \right| #3 \right)}
\newcommand{\rang}{\hbox{rang}}
\newcommand{\rank}{\hbox{rank}}
\newcommand{\signe}{\hbox{signe}}
\newcommand{\sign}{\hbox{sign}}

\newcommand\hA{{\widehat A}}
\newcommand\hB{{\widehat B}}
\newcommand\hC{{\widehat C}}
\newcommand\hD{{\widehat D}}
\newcommand\hE{{\widehat E}}
\newcommand\hF{{\widehat F}}
\newcommand\hG{{\widehat G}}
\newcommand\hH{{\widehat H}}
\newcommand\hI{{\widehat I}}
\newcommand\hJ{{\widehat J}}
\newcommand\hK{{\widehat K}}
\newcommand\hL{{\widehat L}}
\newcommand\hM{{\widehat M}}
\newcommand\hN{{\widehat N}}
\newcommand\hO{{\widehat O}}
\newcommand\hP{{\widehat P}}
\newcommand\hQ{{\widehat Q}}
\newcommand\hR{{\widehat R}}
\newcommand\hS{{\widehat S}}
\newcommand\hTT{{\widehat T}}
\newcommand\hU{{\widehat U}}
\newcommand\hV{{\widehat V}}
\newcommand\hW{{\widehat W}}
\newcommand\hX{{\widehat X}}
\newcommand\hY{{\widehat Y}}
\newcommand\hZ{{\widehat Z}}

\newcommand\ha{{\widehat a}}
\newcommand\hb{{\widehat b}}
\newcommand\hc{{\widehat c}}
\newcommand\hd{{\widehat d}}
\newcommand\he{{\widehat e}}
\newcommand\hf{{\widehat f}}
\newcommand\hg{{\widehat g}}
\newcommand\hh{{\widehat h}}
\newcommand\hi{{\widehat i}}
\newcommand\hj{{\widehat j}}
\newcommand\hk{{\widehat k}}
\newcommand\hl{{\widehat l}}
\newcommand\hm{{\widehat m}}
\newcommand\hn{{\widehat n}}
\newcommand\ho{{\widehat o}}
\newcommand\hp{{\widehat p}}
\newcommand\hq{{\widehat q}}
\newcommand\hr{{\widehat r}}
\newcommand\hs{{\widehat s}}
\newcommand\htt{{\widehat t}}
\newcommand\hu{{\widehat u}}
\newcommand\hv{{\widehat v}}
\newcommand\hw{{\widehat w}}
\newcommand\hx{{\widehat x}}
\newcommand\hy{{\widehat y}}
\newcommand\hz{{\widehat z}}

\newcommand\tA{{\widetilde A}}
\newcommand\tB{{\widetilde B}}
\newcommand\tC{{\widetilde C}}
\newcommand\tD{{\widetilde D}}
\newcommand\tE{{\widetilde E}}
\newcommand\tF{{\widetilde F}}
\newcommand\tG{{\widetilde G}}
\newcommand\tH{{\widetilde H}}
\newcommand\tI{{\widetilde I}}
\newcommand\tJ{{\widetilde J}}
\newcommand\tK{{\widetilde K}}
\newcommand\tL{{\widetilde L}}
\newcommand\tM{{\widetilde M}}
\newcommand\tN{{\widetilde N}}
\newcommand\tOO{{\widetilde O}}
\newcommand\tP{{\widetilde P}}
\newcommand\tQ{{\widetilde Q}}
\newcommand\tR{{\widetilde R}}
\newcommand\tS{{\widetilde S}}
\newcommand\tTT{{\widetilde T}}
\newcommand\tU{{\widetilde U}}
\newcommand\tV{{\widetilde V}}
\newcommand\tW{{\widetilde W}}
\newcommand\tX{{\widetilde X}}
\newcommand\tY{{\widetilde Y}}
\newcommand\tZ{{\widetilde Z}}

\newcommand\ta{{\widetilde a}}
\newcommand\tb{{\widetilde b}}
\newcommand\tc{{\widetilde c}}
\newcommand\td{{\widetilde d}}
\newcommand\te{{\widetilde e}}
\newcommand\tf{{\widetilde f}}
\newcommand\tg{{\widetilde g}}
\newcommand\th{{\widetilde h}}
\newcommand\ti{{\widetilde i}}
\newcommand\tj{{\widetilde j}}
\newcommand\tk{{\widetilde k}}
\newcommand\tl{{\widetilde l}}
\newcommand\tm{{\widetilde m}}
\newcommand\tn{{\widetilde n}}
\newcommand\tio{{\widetilde o}}
\newcommand\tp{{\widetilde p}}
\newcommand\tq{{\widetilde q}}
\newcommand\tr{{\widetilde r}}
\newcommand\ts{{\widetilde s}}
\newcommand\tit{{\widetilde t}}
\newcommand\tu{{\widetilde u}}
\newcommand\tv{{\widetilde v}}
\newcommand\tw{{\widetilde w}}
\newcommand\tx{{\widetilde x}}
\newcommand\ty{{\widetilde y}}
\newcommand\tz{{\widetilde z}}

\newcommand\bA{{\overline A}}
\newcommand\bB{{\overline B}}
\newcommand\bC{{\overline C}}
\newcommand\bD{{\overline D}}
\newcommand\bE{{\overline E}}
\newcommand\bFF{{\overline F}}
\newcommand\bG{{\overline G}}
\newcommand\bH{{\overline H}}
\newcommand\bI{{\overline I}}
\newcommand\bJ{{\overline J}}
\newcommand\bK{{\overline K}}
\newcommand\bL{{\overline L}}
\newcommand\bM{{\overline M}}
\newcommand\bN{{\overline N}}
\newcommand\bO{{\overline O}}
\newcommand\bP{{\overline P}}
\newcommand\bQ{{\overline Q}}
\newcommand\bR{{\overline R}}
\newcommand\bS{{\overline S}}
\newcommand\bT{{\overline T}}
\newcommand\bU{{\overline U}}
\newcommand\bV{{\overline V}}
\newcommand\bW{{\overline W}}
\newcommand\bX{{\overline X}}
\newcommand\bY{{\overline Y}}
\newcommand\bZ{{\overline Z}}

\newcommand\ba{{\overline a}}
\newcommand\bb{{\overline b}}
\newcommand\bc{{\overline c}}
\newcommand\bd{{\overline d}}
\newcommand\be{{\overline e}}
\newcommand\bff{{\overline f}}
\newcommand\bg{{\overline g}}
\newcommand\bh{{\overline h}}
\newcommand\bi{{\overline i}}
\newcommand\bj{{\overline j}}
\newcommand\bk{{\overline k}}
\newcommand\bl{{\overline l}}
\newcommand\bm{{\overline m}}
\newcommand\bn{{\overline n}}
\newcommand\bo{{\overline o}}
\newcommand\bp{{\overline p}}
\newcommand\bq{{\overline q}}
\newcommand\br{{\overline r}}
\newcommand\bs{{\overline s}}
\newcommand\bt{{\overline t}}
\newcommand\bu{{\overline u}}
\newcommand\bv{{\overline v}}
\newcommand\bw{{\overline w}}
\newcommand\bx{{\overline x}}
\newcommand\by{{\overline y}}
\newcommand\bz{{\overline z}}

%
\newcommand{\AAA}{{\cal A}}
\newcommand{\BB}{{\cal B}}
\newcommand{\CC}{{\cal C}}
\newcommand{\DD}{{\cal D}}
\newcommand{\EE}{{\cal E}}
\newcommand{\FF}{{\cal F}}
\newcommand{\GG}{{\cal G}}
\newcommand{\HH}{{\cal H}}
\newcommand{\II}{{\cal I}}
\newcommand{\JJ}{{\cal J}}
\newcommand{\KK}{{\cal K}}
\newcommand{\LL}{{\cal L}}
\newcommand{\NN}{{\cal N}}
\newcommand{\MM}{{\cal M}}
\newcommand{\OO}{{\cal O}}
\newcommand{\PP}{{\cal P}}
\newcommand{\QQ}{{\cal Q}}
\newcommand{\RR}{{\cal R}}
\newcommand{\SS}{{\cal S}}
\newcommand{\TT}{{\cal T}}
\newcommand{\UU}{{\cal U}}
\newcommand{\VV}{{\cal V}}
\newcommand{\WW}{{\cal W}}
\newcommand{\XX}{{\cal X}}
\newcommand{\YY}{{\cal Y}}
\newcommand{\ZZ}{{\cal Z}}
\newcommand{\Id}{operatorname{Id}}

\newcommand{\tbullet}{$\bullet$}
\newcommand{\ot}{\leftarrow}
\newcommand{\carre}{\hfill$\Box$}
\newcommand{\carreb}{\hfill\rule{0.25cm}{0.25cm}}
%
%
\newcommand{\dontforget}[1]
{{\mbox{}\\\noindent\rule{1cm}{2mm}\hfill don't forget : #1
\hfill\rule{1cm}{2mm}}\typeout{---------- don't forget : #1 ------------}}
\newcommand{\note}[2]
{ \noindent{\sf #1 \hfill \today}

\noindent\mbox{}\hrulefill\mbox{}
\begin{quote}\begin{quote}\sf #2\end{quote}\end{quote}
\noindent\mbox{}\hrulefill\mbox{}
\vspace{1cm}
}
\newcommand{\rond}[1]     {\stackrel{\circ}{#1}}
\newcommand{\rondf}       {\stackrel{\circ}{\FF}}
\newcommand{\point}[1]     {\stackrel{\cdot}{#1}}

\newcommand\relatif{{\rm \rlap Z\kern 3pt Z}}

\title{\huge   Pseudoholomorphic discs attached to   $CR$-submanifolds of
  almost complex spaces}
\author{Nikolai Kruzhilin* and Alexandre Sukhov** }
\date{}
\maketitle

{\small {*}Steklov Mathematical Institute, Gubkina Str. 8, GSP-1, 119991,
  Moscow, Russia,  kruzhil@mi.ras.ru

{**}Universit\'e des Sciences et Technologies de Lille, Laboratoire
Paul Painl\'ev\'e,
U.F.R. de
Math\'ematique, 59655 Villeneuve d'Ascq, Cedex, France , 
 sukhov@math.univ-lille1.fr}

\bigskip

Abstract. Let $E$ be a generic real submanifold of an almost complex
manifold. The geometry of  Bishop discs attached to $E$
is studied in terms of the Levi form of $E$.

\bigskip

R\'esum\'e. Nous \'etudions la g\'eom\'etrie des disques de Bishop
attach\'es \`a une sous-vari\'et\'e r\'eelle g\'en\'erique d'une
vari\'et\'e presque complexe.  

\bigskip



Mots-cl\'es: une structure presque complexe, une vari\'et\'e
g\'en\'erique, la forme de Levi, disque de Bishop. 

\bigskip

\section{Introduction}

Attachment of  holomorphic discs to a prescribed
real submanifold of a complex manifold is a well-known and powerfull method of
geometric complex analysis developed by many authors.
Recently, in view of deep connections with  symplectic geometry
discovered by M.Gromov \cite{Gr}, this method found its way to
 the almost complex case. Authors  usually consider the attachment
of pseudoholomorphic discs to totally real submanifolds of almost complex
manifolds, for instance, near points admitting a non-trivial holomorphic tangent
space.

In the present paper  we consider  pseudoholomorphic discs (Bishop
discs) attached to  generic
real submanifolds,  of  a positive complex dimension,
of almost complex manifolds.
  In section 2
we prove the existence of Bishop  discs for  a  generic submanifold
$E$ of
an almost complex manifold $(M,J)$.
 We show that, roughly speaking, these
discs can be parametrized quite similarly to the case of the standard
complex structure. Our proof is based on isotropic dilations of local
coordinates  in a   similar way  to Sikorav's proof of the Nijenhuis-Woolf
theorem on the existence of local pseudoholomorphic discs in a prescribed
direction \cite{Si}.

Our main aim  is to study the geometry of pseudoholomorphic Bishop
discs in terms of the Levi form of a $CR$-submanifold. We begin with
the Levi flat case. In section 3 we prove that each (sufficiently
small) pseudoholomorphic
Bishop disc attached to a real hypersurface with identically
vanishing Levi form  lies  in this hypersurface. Hence
 such a hypersurface contains  pseudoholomorphic discs passing in
an arbitrary prescribed  complex tangent direction (Theorem 3.1).
This gives an
affirmative answer to a   question raised by
Ivashkovich and Rosay \cite{IvRo}.  They
also constructed in  \cite{IvRo} an example of a real hypersurface in an almost
complex
manifold of complex dimension 3 that  has an identically vanishing Levi
form, but contains no complex hypersurfaces. In paricular, this
hypersurface is minimal in the sense of Tumanov \cite{Tu}. Recall that
the well-known result of Tr\'epreau \cite{Tr} and Tumanov \cite{Tu}
claims that in the case of an integrable complex structure, Bishop
discs of a minimal hypersurface fill its one-sided neighborhood. Thus,
Theorem 3.1 in combination  with the  example of Ivashkovich-Rosay shows
that  the Tr\'epreau-Tumanov theorem has no
straightforward generalisation to the almost
complex case.

In section 4  we consider the case  of a
$CR$-submanifold $E$ with  Levi form distinct from zero. We prove that in
this case the corresponding
Bishop discs sweep out a  submanifold  containing   $E$ as an open piece
of the boundary (Theorem 4.1). This  is an almost complex analog
of results due to  Hill--Taiani \cite{HiTa} and Boggess \cite{Bo}, but
our proof in the almost complex setting  requires
 a new idea because
the Nijenhuis tensor (the
torsion) of an almost complex structure has a strong influence on the
geometry of the Levi form of a real submanifold: it is not even always
possible to take a $CR$-submanifold of the standard complex space for a local
model of a $CR$-submanifold of an almost complex space.
To overcome
arising  difficulties, we use in section 4  non-isotropic dilations in a
suitable coordinate system (a similar idea is used in \cite{GoSu} in
order to study boundary behavior of the Kobayashi metric in almost
complex manifolds). It turns out that if  $E$ has
 CR dimension 1, then the non-isotropic dilations allow one to
represent the pair  $(E,J)$ as a small deformation of the pair
$(E_0,J_{st})$, where $E_0$ is the quadric
manifold in $\cc^n$ of which the Levi form with respect to $J_{st}$ coincides
with the Levi form of $E$ with respect to $J$. This results in the
existence of pseudoholomorphic Bishop discs with a certain  special
geometry (Theorem 4.1). The general case  could in principle be treated by
consdering a foliation of $E$ by submanifolds of $CR$-dimenson 1.
However, we  put forward a method allowing us to give a more
straightforward  description  of the  Bishop discs involved in our construction
 in the general case.

We point out that our methods allow one to
deal only with the first Levi form. For instance,  suitable almost complex
 analogs of highly precise results  of Tr\'epreau \cite{Tr} and Tumanov \cite{Tu} require
 another  approach.

This work was partially  carried out  as   the first author
 was visiting  the
University of Lille-1. He thanks this institution for
 hospitality. N.Kruzhilin  is also supported by the RAS Program ''Modern problems of theoretical
 mathematics'' and the Program for Support of Scientific Schools of RF
 (grant N 2040.2003.1).

\section{Existence and local parametrization of Bishop discs}

\subsection{Almost complex manifolds.}
 Let $(M,J)$ be an almost complex manifold
with operator of complex structure $J$.
 Let $\D$ be the unit
disc in $\cc$ and $J_{st}$  the standard (operator of)  complex  structure
on $\cc^n$
for arbitrary  $n$. Let $f$ be
 a smooth map from $D$ into $M$. We say that $f$ is {\it
 $J$-holomorphic}  if $df \circ J' = J \circ df$. We call such a map $f$
a $J$-holomorphic disc and
 denote by   ${\mathcal O}_J(\D,M)$ the set
 of {\it $J$-holomorphic discs} in
 $M$.
We denote by ${\mathcal O}(\D)$ the space of usual holomorphic
 functions on $\D$.

The following lemma shows that an  almost complex manifold
$(M,J)$ can be  locally viewed  as the unit ball $\B$ in
$\cc^n$ equipped with a small almost complex
deformation of $J_{st}$. We shall repeatedly use this
observation in what follows.
\begin{e-lemme}
\label{lemma1}
Let $(M,J)$ be an almost complex manifold. Then for each $p \in
M$,  each  $\delta_0 > 0$, and  each   $k
\geq 0$
 there exist a neighborhood $U$ of $p$ and a
smooth coordinate chart  $z: U \longrightarrow \B$ such that
$z(p) = 0$, $dz(p) \circ J(p) \circ dz^{-1}(0) = J_{st}$,  and the
direct image $z_*(J) := dz \circ J \circ dz^{-1}$ satisfies
the inequality
$\vert\vert z_*(J) - J_{st}
\vert\vert_{\CC^k(\bar {\B})} \leq \delta_0$.
\end{e-lemme}
\proof There exists a diffeomorphism $z$ from a neighborhood $U'$ of
$p \in M$ onto $\B$ such that  $z(p) = 0$ and $dz(p) \circ J(p)
\circ dz^{-1}(0) = J_{st}$. For $\delta > 0$ consider the isotropic dilation
$d_{\delta}: t \mapsto \delta^{-1}t$ in $\cc^n$ and the composite
$z_{\delta} = d_{\delta} \circ z$. Then $\lim_{\delta \rightarrow
0} \vert\vert (z_{\delta})_{*}(J) - J_{st} \vert\vert_{\CC^k(\bar
{\B})} = 0$. Setting $U = z^{-1}_{\delta}(\B)$ for sufficiently small positive
$\delta$  we obtain the required result.

{\bf The operators $\partial_J$ and $\bar{\partial}_J$ .}

Let $(M,J)$ be an almost complex manifold. We denote by $TM$ the real
tangent bundle of $M$ and by $T_{\cc} M$ its complexification. Recall
that $T_{\cc} M = T^{(1,0)}M \oplus T^{(0,1)}M$ where
$T^{(1,0)}M:=\{ X \in T_{\cc} M : JX=iX\} = \{\zeta -iJ \zeta, \zeta \in
TM\},$
and $T^{(0,1)}M:=\{ X \in T_{\cc} M : JX=-iX\} = \{\zeta +
iJ \zeta, \zeta \in TM\}$.
 Let $T^*M$ be  the cotangent bundle of  $M$.
Identifying $\cc \otimes T^*M$ with
$T_{\cc}^*M:=Hom(T_{\cc} M,\cc)$ we define the set of complex
forms of type $(1,0)$ on $M$ as
$
T_{(1,0)}M=\{w \in T_{\cc}^* M : w(X) = 0, \forall X \in T^{(0,1)}M\}
$
and we denote  the set of complex forms of type $(0,1)$ on $M$ by
$
T_{(0,1)}M=\{w \in T_{\cc}^* M : w(X) = 0, \forall X \in T^{(1,0)}M\}
$.
Then $T_{\cc}^*M=T_{(1,0)}M \oplus T_{(0,1)}M$.
This allows us to define the operators $\partial_J$ and
$\bar{\partial}_J$ on the space of smooth functions  on
$M$~: for a  smooth complex function $u$ on $M$ we set $\partial_J u =
du_{(1,0)} \in T_{(1,0)}M$ and $\bar{\partial}_Ju = du_{(0,1)}
\in T_{(0,1)}M$. As usual,
differential forms of any bidegree $(p,q)$ on $(M,J)$ are defined
by  exterior multiplication.

{\bf Plurisubharmonic functions.} We say that
an upper semicontinuous function $u$ on $(M,J)$ is
{\it $J$-plurisubharmonic} on $M$ if the composition $u \circ f$
is subharmonic on $\Delta$ for every $f \in {\mathcal O}_J(\D,M)$.

 Let $u$ be a $\CC^2$ function on $M$, let $p \in M$
and $v \in T_pM$. {\it The Levi
form} of $u$ at $p$ evaluated on $v$ is defined by the equality
$L^J(u)(p)(v):=-d(J^* du)(X,JX)(p)$
where $X$ is an arbitrary vector field on $TM$ such that $X(p) = v$ (of course,
this definition is independent of one's choice of $X$).

 The following result is well known (see, for instance, \cite{IvRo}).
\begin{e-pro}\label{PROP1}
Let $u$ be a $\CC^2$ real valued function on $M$, let $p \in M$ and $v \in
T_pM$. Then $ L^J(u)(p)(v) = \Delta(u \circ f)(0)$
where $f$ is an arbitrary $J$-holomorphic disc in $M$ such that
 $f(0) = p$ and $
df(0)(\partial / \partial Re \zeta) = v$ (here $\zeta$ is      the
standard complex coordinate variable  in $\cc$).
\end{e-pro}

The Levi form is obviously invariant with respect to biholomorphisms.
More precisely, let $u$ be a $\CC^2$ real valued function on $M$,
let $p \in M$ and $v \in T_pM$.
If $\Phi$ is a
$(J,J')$-holomorphic diffeomorphism from $(M,J)$ into $(M',J')$, then
$ L^J(u)(p)(v) =
 L^{J'}(u \circ \Phi^{-1})(\Phi(p))(d\Phi(p)(v))$.

Finally, it follows from Proposition~\ref{PROP1} that a
$C^2$-smooth real function $u$ is
$J$-plurisubharmonic on $M$ if and only if $ L^J(u)(p)(v) \geq 0$
for all   $p \in M$, $v \in T_pM$.
Thus, similarly to  the case of
an integrable structure one arrives in a natural way to
the following definition: a $\CC^2$ real valued
function $u$ on $M$ is {\it strictly
$J$-plurisubharmonic} on $M$ if  $ L^J(u)(p)(v)$
is positive for every $p \in M$, $v \in T_pM \backslash \{0\}$.

It follows easily from Lemma \ref{lemma1} that
for every point $p\in (M,J)$ there exists a neighborhood $U$ of $p$
 and a diffeomorphism $z:U \rightarrow \B$ with
center   at $p$ (in the sense that $z(p) =0$) such that the function
$|z|^2$ is
$J$-plurisubharmonic on $U$ and $z_*(J) = J_{st} + O(\vert z \vert)$.

Let $u$ be a $\CC^2$ function   in a neighborhood
of a point $p$ of $(M,J)$  that is  strictly $J$-plurisubharmonic.
Then there exists a neighborhood $U$ of $p$  with local
complex coordinates $z:U \rightarrow \B$  such
that the function $u - c|z|^2$ is $J$-plurisubharmonic on $U$ for some
constant $c > 0$.

{\bf Real submanifolds in almost complex manifolds.} Let $E$ be a real
submanifold of codimension $m$
in an almost complex manifold $(M,J)$ of complex dimension $n$. For every $p$ we
denote by $H_p^J(E)$ the maximal complex (with respect to $J(p)$)
subspace of the tangent space $T_p(E)$. Similarly to the integrable
case, $E$ is said to be a CR manifold if the (complex) dimension of
$H_p^J(E)$ is independent on $p$; it is called the CR dimension of $E$
and is denoted by $CR dim E$.

In complex analysis by a {\it generic\/}
submanifold of a complex manifold   one usually means
a submanifold $E$ such that                at
every point $p \in E$ the complex linear span of $T_p(E)$ coincides with
the tangent space of the ambient manifold. We think that the use of this
term in precisely that sense  outside the framwork of complex analysis
proper can sometimes be misleading.
For this reason we shall provisionally call submanifolds with similar
properties of an almost complex manifold $M$ (that is,
 submanifolds such that  the complex linear span of $T_p(E)$ at each point
 coincides with $TM$) {\it generating\/} submanifolds.
Of course,   every generating
submanifold is CR.

If $E$ is defined as the common zero level of functions $r_1,...,r_m$,
then after the standard identification of $TM$ and $T^{(1,0)}M$
 $H_p^J(E)$ can be defined as  the zero subspace of the forms
${\partial}_Jr_1$, ... ${\partial}_Jr_m$.
In particular, let $E$ be a smooth real  hypersurface in an almost
complex manifold $(M,J)$ defined by an equation $r = 0$.
We say that $E$ is
strictly pseudoconvex (for $J$) if the Levi form of $r$ is strictly positive
definite on each holomorphic tangent space $H_p^J(E)$, $p \in
\Gamma$. Of course, this definition does not depend on one's choice of the
defining function $r$. We shall require the following result,
which is  well-known
in the case of an integrable structure.

\begin{e-lemme}
\label{strictlypsh}
Let $E$ be a strictly pseudoconvex hypersurface in an almost
complex manifold $(M,J)$. Then $E$ admits a strictly
plurisubharmonic defining function in a neighborhood of each its point
$p$.
\end{e-lemme}

The proof is quite similar to the case of the standard structure. Selecting
suitable  local coordinates $Z = (z,w_1,...,w_{n-1})$ in $\cc^n$ we can assume that $M$ is a neighborhood of the
origin in $\cc^n$ and $p = 0$; as usual, we assume that $J(0) =
J_{st}$. Furthermore, we may suppose that $r(Z) = z + \overline z
+ O(\vert Z \vert^2)$ and therefore $H_0^J(E) = \{ z = 0
\}$. Then $L^J_0(r^2)(v) = \vert v_1 \vert^2$ for  complex tangent
vectors $v \in T_0(M)$. Since $L_0^J(r)$ is strictly positive definite
on $H_0^J(E)$, the Levi form $L_0^J(r + Cr^2)$ is strictly
positive on $T_0(M)$ for a sufficiently large positive  constant $C$.

\subsection{Bishop discs and Bishop's equation} Let $(M,J)$ be a
smooth almost complex
manifold of real dimension
$2n$ and $E$  a generating submanifold of $M$ of real codimension
$m$. A $J$-holomorphic disc $f:\D \longrightarrow M$ continuous on
$\overline \D$ is called a {\it Bishop disc} if $f(b\D) \subset E$
(where $b\D$ denotes the boundary of $\D$). Our aim  is to prove the
existence and to describe certain classses of Bishop discs attached to $E$.

Consider the case where $E$ is defined as the zero set of an
$\R^m$-valued  function $r = (r^1,...,r^m)$  on $M$.
Then a smooth map $f$
defined on $\D$ and continuous on
$\overline\D$ is a Bishop disc if and only if it satisfies the
following non-linear boundary problem of the
Riemann-Hilbert type for the
quasi-linear operator $\overline\partial_J$:

\begin{displaymath}
(RH): \left\{ \begin{array}{ll}
\overline\partial_J f(\zeta) = 0, \zeta \in \D\\
r(f)(\zeta) = 0, \zeta \in b\D
\end{array} \right.
\end{displaymath}

To describe  solutions of this problem we
fix a chart $U \subset M$ and a coordinate diffeomorphism $z: U
\longrightarrow \B^n$ where $\B^n$ is  the unit ball
of $\cc^n$. Identifying $M$ with $\B^n$
we may assume that in these coordinates $J= J_{st} + O(\vert z \vert)$
and the norm $\parallel J - J_{st} \parallel_{C^k(\overline \B^n)}$ is
small enough for some  positive real
$k$ in accordance with Lemma \ref{lemma1}. (Here $k$ can be arbitrary,
but we  assume it for convenience to be non-integer
and fix it throughout what follows.)
 More precisely,
using the notation $Z = (z,w)$, $z = (z_1,...,z_m)$,
$w = (w_{1},...,w_{n-m})$ for the standard coordinates in $\cc^n$, we may also
assume that $E \cap U$ is described by the equations

\begin{eqnarray}
\label{manifold}
r(Z) = Re z - h(Im z,w) = 0
\end{eqnarray}
with  vector-valued $C^{\infty}$-function $h:\B \longrightarrow \R^m$ such that
$h(0) = 0$ and $\bigtriangledown h(0) = 0$.

Similarly to  the proof of Lemma \ref{lemma1} consider the isotropic dilations
$d_{\delta}: Z \mapsto Z' = \delta^{-1}Z$. In the new
$Z$-variables  (we drop the primes) the image
$E_{\delta} = d_{\delta}(E)$ is
defined by the equation $r_\delta (Z):=
\delta^{-1}r(\delta Z) = 0$. Since the function
$r_\delta$ approaches $Re z$ as $\delta \longrightarrow 0$,
the manifolds $E_\delta$ approach the flat manifold $E_0=\{
Re z = 0\}$, which, of course, may be identified with the real tangent
space to $E$ at the origin. Furthermore, as         seen in the proof
of Lemma \ref{lemma1}, the structures $J_\delta :=
(d_\delta)_*(J)$
converge to $J_{st}$ in the $C^k$-norm  as
$\delta \longrightarrow 0$. This allows us to find explicitly   the
$\overline\partial_J$-operator in the $Z$ variables.

Consider a $J_\delta$-holomorphic disc
$f: \D \longrightarrow (\B^n,J_\delta)$. The
$J_\delta$-holomorphy condition $J_\delta(f) \circ f_* = f_* \circ J_{st}$ can be
written in the following form.

\begin{eqnarray}
\label{CR}
\frac{\partial f}{\partial \overline\zeta} + Q_{J,\delta}(f)
\left ( \frac{\partial \overline f}{\partial \overline \zeta} \right )
= 0
\end{eqnarray}
where $Q_{J,\delta(Z)}$ is the
complex $n\times n$ matrix  of an operator the composite of which with complex
conjugation is equal to the endomorphism
$- (J_{st} + J_\delta(Z))^{-1}(J_{st} - J_\delta(Z))$
(which is an anti-linear operator with respect to the standrard
structure $J_{st}$). Hence the entries of the matrix $Q_{J,\delta}(z)$
are smooth functions of $\delta,z$ vanishing identically in $z$
for $\delta = 0$.

Using the Cauchy-Green transform

$$T_{CG}(g) = \frac{1}{2\pi i} \int\int_{\D} \frac{g(\tau)}{\zeta -
  \tau}d\tau \wedge d\overline\tau$$
we may write   $\overline\partial_J$-equation (\ref{CR}) as follows:

$$\frac{\partial}{\partial\overline \zeta} \left ( f + T_{CG}\left
(Q_{J,\delta}(f)
\left ( \frac{\partial \overline f}{\partial \overline \zeta} \right
)\right ) \right ) = 0$$

 According to  classsical results \cite{Ve}, the Cauchy-Green
 transform is a continuous linear operator from $C^k(\overline \D)$ into
 $C^{k+1}(\overline \D)$ (recall that $k$ is noninteger). Hence  the operator

$$\Phi_{J,\delta}: f \longrightarrow g =  f + T_{CG}\left (Q_{J,\delta}(f)
\left ( \frac{\partial \overline f}{\partial \overline \zeta} \right
)\right )$$
takes the space   $C^{k}(\overline\D)$  into itself.
Thus, $f$ is $J_\delta$-holomorphic if
and only if $\Phi_{J,\delta}(f)$ is holomorphic (in the usual sense)
on $\D$.
For sufficiently small positive $\delta$  this is an
invertible operator on a neihbbourhood  of zero in $C^k(\overline D)$
which  establishes a one-to-one correspondence between the
sets of $J_\delta$-holomorphic and holomorphic discs  in $\B^n$.

These considerations allow us to replace the non-linear Riemann-Hilbert problem
(RH) by {\it generalized Bishop's equation}

\begin{eqnarray}
\label{Bishop2}
r_\delta(\Phi_{J,\delta}^{-1}(g))(\zeta) = 0, \zeta \in b\D
\end{eqnarray}
for  an unknown {\it holomorphic} function $g$ in the disc
 (with respect to the standard complex structure).

If $g$ is a solution of the boundary  problem (\ref{Bishop2}), then $f =
\Phi_{J,\delta}^{-1}(g)$ is a Bishop disc with  boundary
attached to $E_\delta$. Since the manifold $E_\delta$ is
biholomorhic via  isotropic dilations to the initial manifold $E$,
the solutions of the equation (\ref{Bishop2}) allow to describe
Bishop's discs attached to $E$. Of course, this gives just the discs
close enough (in the $C^k$-norm) to the trivial solution $f \equiv 0$
of the problem (RH).

\subsection{Solution of  generalized Bishop's equation}

Let  ${\cal U}$ be a  neighborhood of the origin in $\R$,
$X'$ a sufficiently small neighborhood of the origin in the Banach  space $ ({\cal O}(\D) \cap C^{k}(\overline \D))^m$  (with
 positive noninteger $k$), $X''$  a neighborhood of the origin in the Banach space  $ ({\cal
  O}(\D) \cap C^k(\overline \D))^{n-m}$, and $Y$ the Banach space
$(C^{k}(b\D))^m$. If $z \in X'$, $z:\zeta \mapsto z(\zeta)$ and
$w \in X''$, $w:\zeta \mapsto w(\zeta)$ are holomorphic
discs, then we denote by $Z$ the holomorphic disc $Z = (z,w)$.
We may also assume that $J$ is a $C^{2k}$-smooth
real $(2n \times 2n)$-matrix valued function
 and denote by $W$  the Banach space of these functions.

Consider the map of Banach spaces
$R:  W \times X'\times X'' \times {\cal U}  \longrightarrow Y$
defined as follows:

$$R : (J,z,w,\delta) \mapsto
r_\delta(\Phi_{J,\delta}^{-1}(z,w))(\bullet)\vert b\D.$$

Let $\phi$ be a   $C^{2k}$-map between two domains in $\R^n$ and $\R^m$ ; it determines a map
$\omega_\phi$ acting   by compostion on  $C^k$-smooth maps $g$ into
the source domain:
$\omega_\phi: g \mapsto \phi(g)$. The well-known fact is that
$\omega_\phi$ is a $C^k$-smooth map between the corresponding spaces
of $C^k$-maps. In our case this means that the map $R$ is of class $C^k$.
For a holomorphic disc $Z = (z,w)  \in X' \times X''$ and $J \in W$
the tangent map $D_{X'}R(J,z,w,0): X'
\longrightarrow Y$ (the partial derivative with respect
to the space $X'$) is defined by the equality
 $D_{X'}R(J,z,w,0)(q) = (Re q_1,...,Re q_m)$.
Clearly, it is surjective and has the kernel of real dimension $m$
consisting of constant functions $h = i(c_1,...,c_m)$, $c_j \in \R$.
(We point out that the map $D_{X'}R$  coincides with the
tangent map arising after the linearization of Bishop's equation
corresponding to the standrard structure $J_{st}$.) Therefore, by the
implicit function theorem \cite{Ze}  there exists  $\delta_0 > 0$, a
neighborhood $V_1$ of the origin in $X'$, a neighborhood $V_2$ of the
origin in $X''$, a neighborhood $V_3$ of the origin in $\R^m$, a
neigborhood $W_1$ of $J_{st}$ in $W$  and a $C^k$
smooth  map $G: W_1 \times V_2 \times V_3 \times [0,\delta_0] \longrightarrow
V_1$ such
that for every $(J,w,c,\delta) \in W_1 \times V_2 \times V_3 \times
[0,\delta_0]$ the function $g = (G(J,w,c,\delta)(\bullet),w(\bullet))$ is
the unique solution of  generalized Bishop's equation (\ref{Bishop2})
belonging to $V_1 \times V_2$.

Now, the pullback $f = \Phi_{J,\delta}^{-1}(g)$ gives us
a $J_\delta$-holomorphic disc attached
to $E_\delta$.   Thus, the  initial data consisting of $J \in W_1$,
and a set $(c_1,...,c_m,w)$, $c_j \in \R$, $w \in V_2$ define
for each small $\delta$ a unique
$J_\delta$-holomorphic disc $f$
attached to $E_\delta$. Since the almost complex structures $J$
and $J_\delta$ are biholomorphic via  isotropic dilations, we can
give   the following description of local solutions of  Bishop's
equation (the problem (RH)).

\begin{e-theo}
\label{theosolvingBishop}  Let $E$ be a smooth submanifold of $\cc^n$
defined as the zero set of a smooth $\R^m$-valued function $r$ of the
form     (\ref{manifold}). Then
there exists  a
neighborhood $U$  of  the origin in $(C^k(\overline\D))^n$,    a
neighborhood $W_1$ of $J_{st}$ in the space $W$,
a neighborhood $V_2$ of the
origin in $X''$,
and  a neighborhood $V_3$ of the origin in $\R^m$
 such that for each $J\in  W_1$
the set of maps  in  $U$ that are Bishop discs
attached to $E$ with respect to $J$
  is a   Banach submanifold of class $C^k$ in  $U$ with local chart
defined by a (smooth)  map $F:  V_2 \times V_3  \longrightarrow
U$, which  depends smoothly  on $J$.
\end{e-theo}

The proof follows from the above analysis of the Bishop equation. One
merely fixes    some value of $\delta$, $0<\delta\leq \delta_0$
and observes that the families
of  Bishop discs corresponding to distinct values of
$\delta\neq 0$ are taken into one another by the corresponding dilations.

One important consequence of this statement is as follows: if $E_1$ and $E_2$
are $C^{2k}$-close submanifolds of $\cc^n$ defined by equations of the form
(\ref{manifold}) and $J_1$ and $J_2$ are $C^{2k}$
close almost complex structures, then  there exists a (locally defined)
diffeomorphism between the corresponding    $C^{k}$ Banach submanifolds
of Bishop discs that depends smoothly on the pairs $E_j$, $J_j$, $j=1,2$,
and is the identity  in the case of equal pairs.

\section{Hypersurfaces with vanishing Levi form}

Here we prove the following result.

\begin{e-theo}
Let $\Gamma$ be a real hypersurface in an almost complex manifold
$(M,J)$ with  Levi form vanishing identically at the   points of
$\Gamma$. Then at each  point $p \in \Gamma$ and for each  direction $v \in
H_p(\Gamma)$ there exists a $J$-holomorphic disc $f:\D \longrightarrow
M$ such that $f(0) = p$, $df(0)(\partial/\partial Re \zeta) = v$ and
$f(\D) \subset \Gamma$.
\end{e-theo}

If $M$ has real dimension 4, then this result can be proved in the same
fashion as in complex analysis, by the application of Frobenius's theorem
to complex tangent spaces to $\Gamma$. However, this does not work in higher
dimensions, when in the case of an almost complex structure distinct from
the standard one the distribition of the planes $H_p(\Gamma)$ is not
necessarily involutive.
The idea of our proof is to show that   {\it each} (sufficiently small) Bishop
 disc for $\Gamma$ lies
in $\Gamma$.

As before, passing to local coordinates $Z = (z,w_1,...,w_{n-1})$
we may assume that $\Gamma$ is a real hypersurface in a
neighborhood $\Omega$ of the origin in $\cc^n$ and $J(0) = J_{st}$.
Let $r$ be a local defining function of
$\Gamma$ in  $\Omega$. Denote by $\Omega^+$ (resp. $\Omega^-$) the
 domain $\{ Z \in \Omega: r(Z) > 0 \}$ (resp. $\{ Z \in \Omega: r(Z) <
 0 \}$).

A
neighborhood $\Omega$  is supposed  to be small enough; in particular,
we may assume that

(a) the function $\vert Z \vert^2$ is strictly
$J$-plurisubharmonic on $\Omega$ and there exists a constant
$\varepsilon_0 > 0$ such that the
value of the Levi form (with respect to $J$) of the function
$\vert Z \vert^2$ on a vector $v$ at a point
$p \in \Omega$ is minorated by $(\varepsilon_0/2)\parallel v \parallel^2$.

 For a constant $N >1$, which will be chosen later, and
sufficiently small $\varepsilon > 0$
 consider the function
$r_{\varepsilon}(Z) = r(Z) + \varepsilon \vert Z \vert^2 -
\varepsilon/N$ and the hypersurface $\Gamma_{\varepsilon}:= \{ Z \in
\Omega :
r_{\varepsilon}(Z) = 0 \}$.  Recall that $\B^n$ is     the unit ball
in $\cc^n$; we may assume that

(b) the ball $\B^n$ lies
 in $\Omega$.

We shall make our choice of a coordinate system
 more  precise. Namely, after a
 $\cc$-linear change of coordinates preserving the previous
 assumptions, we may assume that

(c) $r(Z) = x -h(y,w)$ on $\Omega$ (as usual, $z = x+iy$ and $w =
 (w_1,...,w_{n-1}$)); in particular  $H_0^J(\Gamma) = \{Z: z = 0 \}$.

and

(d)  for every $Z \in \Omega$ and for
$\varepsilon < 1$
   the kernel of the form $\partial_J
r_\varepsilon(Z)$ is in a one-to-one correspondence
with $H_0^J(\Gamma)$ via the
projection $(z,w) \mapsto (0,w)$.

Furthemore, there exists a constant
$c > 0$ such that for every $z \in \Omega$
one has  $c^{-1}dist(Z,\Gamma) \leq \vert r(Z) \vert \leq c\,
dist(Z,\Gamma)$ (where $dist$ is      the Euclidean distance).
For $Z \in (1/2N)\B^n \cap \Gamma_{\varepsilon}$
we have $\vert r(Z) \vert = \varepsilon/N - \varepsilon \vert Z
\vert^2$ so that
\begin{eqnarray}
\label{*}
\varepsilon/2N \leq  \vert r(Z) \vert \leq \varepsilon/N
\end{eqnarray}

Throughout the rest of the proof we shall stay in $ (1/2N)\B^n$.
Note that

(e)
pieces of            the
hypersurfaces $\Gamma_{\varepsilon}$ and $\Gamma_{-\varepsilon}$ form
 a foliation of
 $ (1/2N)\B^n$.

\begin{e-lemme}
If   $N$ is a sufficiently  large fixed  constant, then for sufficiently
small
 $\varepsilon > 0$ the hypersurface
$\Gamma_{\varepsilon}$ is strictly $J$-pseudoconvex at all its point lying
 in the ball  $(1/4N)\B^n$.
\end{e-lemme}
\proof   For every $Z \in \Omega^+ \cap (1/2N)\B^n $ there exists
unique $\varepsilon = \varepsilon(Z)$ such that $Z \in
\Gamma_{\varepsilon}$; clearly, the function $Z \mapsto
\varepsilon(Z)$ is smooth on $(1/2N)\B^n$. Therefore
for every $Z \in
\Gamma_{\varepsilon} \cap (1/2N)\B^n$
the value of the Levi form $L_Z^J(r)(v)$ of
$r$ at $Z$ at a vector $v \in
H_Z^J(\Gamma_{\varepsilon})$ has the estimate
\begin{eqnarray}
\label{est}
\vert L_Z^J(r)(v) \vert < c_1(d_1+d_2)\parallel v \parallel^2,
\end{eqnarray}
where $d_1$ is the distance from $Z$ to the closest point
$Z_0$ on the hypersurface $\Gamma = \Gamma_0$
and $d_2$ is the distance between  $H_Z^J(\Gamma_{\varepsilon})$ and
$H_{Z_0}^J(\Gamma_{0})$ measured in some smooth metric on the corresponding
Grassmanian. Indeed, since $\vert L_Z^J(r)(v) \vert = \vert
L_Z^J(r)(v/\parallel v \parallel) \vert \parallel v \parallel^2$, it   is
sufficient to find an estimate of $ \sup \{ \vert L_Z^J(r)(u)
\vert : u \in H_Z^J(\Gamma_\varepsilon), \parallel u \parallel =1
\}$. Consider the real unit spheres $S(Z) = \{ u \in
H_Z^J(\Gamma_\varepsilon): \parallel u \parallel =1\}$ and $S(Z_0) =
\{ u' \in
H_{Z_0}^J(\Gamma): \parallel u' \parallel =1\}$ in the tangent spaces
$H_Z^J(\Gamma_\varepsilon)$ and $H_{Z_0}^J(\Gamma)$ respectively.
In what follows the Levi forms
$L^J_Z(r)$ and $L^J_{Z_0}(r)$ are viewed  as quadratic forms on
$\R^{2n}$ since the local coordinates are fixed; the tangent spaces
$H_Z^J(\Gamma_\varepsilon)$ and $H_{Z_0}^J(\Gamma)$ are identified with
subspaces in $\R^{2n}$. Denote by $\hat L^J_Z(r)$ the polarization of
$L^J_Z(r)$, that is, the corresponding  bilinear form on
$\R^{2n}$.

 For any vector $u \in S(Z)$  we have
$$\vert L_Z^J(r)(u) \vert \leq \inf_{ u' \in S(Z_0)} \left (\vert
L_Z^J(r)(u') \vert +  \vert L_Z^J(r)(u - u')\vert  + 2\vert \hat
L_Z^J(r)(u',u-u') \vert \right ). $$   When $Z \in \Gamma$
the form $L_Z^J(r)(u')$ vanishes for any $u'\in H_Z^J(\Gamma)$; so there exists a
constant $c_1$  such that $$\sup
\{ \vert
L_Z^J(r)(u') \vert: u' \in S(Z_0) \} \leq C d_1.$$
Furthemore, there exist constants $c'_1$ and $c'_2$ such that
$$\sup_{u \in S(Z)}\left ( \inf_{ u' \in S(Z_0)}
 \vert L_Z^J(r)(u - u')\vert\right )  \leq c'_1 \sup_{u \in S(Z)}
\left ( \inf_{ u' \in S(Z_0)}  \parallel u - u' \parallel^2 \right )
 \leq c'_2d_2^2.$$
In a similar way

$$\sup_{u \in S(Z)}\left ( \inf_{ u' \in S(Z_0)} 2\vert \hat
L_Z^J(r)(u',u-u') \vert \right ) \leq  c''_1 \sup_{u \in S(Z)}
\left ( \inf_{ u' \in S(Z_0)}  \parallel u - u' \parallel \right )
 \leq  c''_2d_2$$
for some positive constants $ c''_1$ and $ c''_2$.
 Obviously, $d_1\sim \vert r \vert$. Moreover, a direct estimate from above of the
 quantities  $\vert \partial r_\varepsilon/\partial z_k(Z) - \partial
r/\partial z_k(Z_0)\vert$ and $\vert\partial r_\varepsilon/\partial \overline z_k(Z) - \partial
r/\partial \overline z_k(Z_0)\vert$ shows that  $d_2<c_2(d_1
+\varepsilon |Z|)$.

Observing  that if $Z\in\Gamma_\varepsilon$, then
$\varepsilon=r(Z)/(\frac{1}{N}-|Z|^2)$,
and also that $|Z|<1/2N$ we see that

\begin{eqnarray}
\label{**}
\vert L_Z^J(r)(v) \vert \leq  (c_3 \varepsilon(Z)/N)\parallel v \parallel^2
\end{eqnarray}
We point out that the constant $c_3$ is independent of
$\varepsilon$ and $N$. Fix $N \geq \max
\{2, 4
c_3/\varepsilon_0 \}$. Then,  in view of  condition (a) and (\ref{**}),
  the Levi form of $r_{\varepsilon}$ is strictly positive on
$H_Z^J(\Gamma_{\varepsilon})$ which proves the lemma.

By  Theorem \ref{theosolvingBishop}, in each  sufficiently small
 neighborhood of
the origin  there exists
a family of $J$-holomorphic Bishop discs with boundaries in
$\Gamma$. Fix a such a neighborhood $U \subset (1/2N\B^n)$.

\begin{e-lemme}
\label{embedding}
Let $f: \D \longrightarrow U$ be a $J$-holomorphic Bishop
disc, that is,  let  $f$ be a pseudoholomorphic map
 continuous on $\overline \D$ such that $f(b\D) \subset
\Gamma$. Then $f(\overline\D)$ lies in $\Gamma$.
\end{e-lemme}

\proof Assume by contradiction that $f(\D)$  does not lie  in
$\Gamma$. Recall  that $U^+:= U \cap \{ r > 0 \}$ is filled by strictly
pseudoconvex hypersurfaces $\Gamma_{\varepsilon}$, $0 \leq \varepsilon
\leq \varepsilon_0$.
We may assume that there exists a connected  open subset $G $ of $  \D$ such
that $f(G)$ lies  in
$U^+$ (otherwise we replace $r$ by $-r$)
and $f(bG) \subset \Gamma$. Consider the set $A = \{ \varepsilon > 0: (r_\varepsilon \circ
f) \vert_G < 0 \}$. This set is not empty if the disc $f$ is small
enough. Let $\varepsilon_1 = \inf A$. Then $\varepsilon_1 > 0$, the
hypersurface $\Gamma_{\varepsilon_1}$ is strictly pseudoconvex,
 $(r_{\varepsilon_1} \circ f) \vert_{bG} < 0$ and $r_{\varepsilon_1}\circ
f\vert_G \leq 0$. Moreover, there exists an interior point $\zeta \in
G$
such that $r_{\varepsilon_1} \circ f (\zeta) = 0$.

On the other hand
 by Lemma \ref{strictlypsh} the hypersurface $\Gamma_{\varepsilon_1}$ admits a
strictly plurisubharmonic defining function in a neighborhood of 
the point
$f(\zeta)$. This contradicts    the maximum principle and proves the
lemma.

We now can     prove Theorem 3.1. Similarly to the previous
sections consider the isotropic dilations $d_\delta: Z \mapsto Z'=
\delta^{-1}Z$. The image $\Gamma_\delta: = (d_\delta)_*(\Gamma)$ of
  the hypersurface $\Gamma$ approaches the hyperplane $\Gamma_0 = \{ Re
  z = 0 \}$ as $\delta \longrightarrow 0$. Let $U_j$, $j=1,2$ be
  neighborhoods of the origin in $\cc^{n-1}$ and $U_3$ a
  neighborhood of the origin in $\R$; we assume that
 these neighborhoods are
  {\it sufficiently small}. For $p \in U_1$, $v \in U_2$ and
  $c \in U_3$ consider a $J_{st}$-holomorphic disc
$f(p,v,c)(\zeta) = (ic, p + v\zeta)$ that is a Bishop disc lying  in the
  hyperplane $\Gamma_0$. The   centers $f(p,v,c)(0)$ of such discs fill a
  neighborhood of the origin in $\Gamma_0$ and their tangent vectors
  (at centers)
$df(p,v,c)(0)(\partial/\partial Re \zeta)$ fill a neighborhood of the
  origin in the holomorphic tangent space $H_q(\Gamma_0)$ for any $q
  \in \Gamma_0$ in a neighborhhod of the origin.
By Theorem \ref{theosolvingBishop} for any $\delta > 0$ there exists a
family of discs $F(\delta,p,v,c)(\bullet)$ smoothly depending on
parameters $\delta,p,v,c)$ such that

\begin{itemize}
\item[(a)] every disc
$F(\delta,p,v,c)(\bullet)$ is $J_\delta$ holomorphic (where as usual
$J_\delta$ denotes the direct image $(d_\delta)_*(J))$;
\item[(b)] for every sufficiently small
positive $\delta$              every disc
  $F(\delta,p,v,c)(\bullet)$ is a Bishop disc for $\Gamma_\delta$,
 that is,   $F(\delta,p,v,c)(b\D) \subset \Gamma_\delta$;
\item[(c)] we have $F(0,p,v,c)(\bullet) = f(p,v,c)(\bullet)$,
 so that   the family $\{ F(\delta,p,v,c)(\bullet) \}$ of
  $J_\delta$-holomorphic discs is a small deformation of the family
  $\{ f(p,v,c)(\bullet) \}$.
\end{itemize}

By Lemma \ref{embedding}, for small $\delta > 0$
every disc $F(\delta,p,v,c)(\D)$ lies in $\Gamma_\delta$. By standard
arguments  their centers fill a neighborhood $U$ of the origin on
$\Gamma_\delta$ and at every point $z \in  U$  their tangent
vectors fill a neighborhood of the origin in the tangent space
$H_z^{J_\delta}(\Gamma_\delta)$. Since the structures $J_\delta$ and
$J$ are biholomorphic, the  proof of the theorem is complete.

\section{Manifolds with  non-trivial  Levi form}

In this section we prove the following  result

\begin{e-theo}
\label{theo1}
Let $E = \{ r: = (r^1,...,r^m) = 0, j=1,...,m \}$ be a (germ of a) smooth
generating submanifold
passing through a point $p$ in an almost complex
manifold $(M,J)$. Suppose that there exists $j$ and a vector $v \in
H_p^J(E)$ such that the Levi form  $L^J_p(r^j)(v)$  does not
vanish. Then for fixed
 non-integer $k > 2$ there exists
 in a neighborhood of $p$ a $C^{k}$ smooth generating
 manifold $\tilde E$ of dimension
$\dim E + 1$ with 
boundary such that every point of $\tilde E$
belongs to a
  $J$-holomorphic disc with boundary on $E$ and  $E$ is
   the boundary of $\tilde E$.
\end{e-theo}

Our proof is based on  non-isotropic scaling.     Isotropic
dilations used in the previous section can not be applied here since
they do not give  one  the control over the Levi form of $E$. The
crucial technical point here is a choice of a suitable coordinate
system ``normalizing'' an almost complex structure.
 Indeed, the following elementary example shows
the basic difficulty in dealing with the almost complex case if a
coordinate system is not good enough. Consider
in $\cc^2$ the real hyperplane $\Pi: Re z_2 = 0$, which is Levi flat
in  the standard complex structure $J_{st}$ of $\cc^2$.
(Throughout, we identify an almost complex structure on a manifold with the corresponding
field of operators on the tangent space.)
Consider the
diffeomorphism $\Phi: (z_1,z_2) \mapsto (z_1, z_2 - \vert z_1
\vert^2)$. The image $\Phi(\Pi)$ is the hypersurface $\Gamma: Re z_2 +
\vert z_1 \vert^2 = 0$ and the direct image of the standrard structure
is the almost complex structure $J(\Phi(z)) = d \Phi(z) \circ J(z) \circ
d \Phi^{-1}(z)$. The structure $J$ coincides with $J_{st}$ at the origin,
so that  $J(z) = J_{st} + O(\vert z \vert)$ and  the hypersurface $\Gamma$ is
strictly pseudoconvex with respect to $J_{st}$, but Levi flat with
respect to $J$\,!

\subsection{ The case  $CRdim E = 1$}

We begin with this case since it  is particularily convenient for
non-isotropic dilations. Passing to suitable local coordinates
(similarly to  the previous section
we use the notation $Z = (z_1,...,z_{n-1},w)$) we may
assume that $M$ is a neighborhood of the origin in $\cc^n$
and $J$ is a
smooth matrix valued function of the form $J = J_{st} + O(\vert Z
\vert)$.
 Moreover, we may assume that the holomorphic tangent space $H_0^J(E)$
 coincides with
the line $l= (0,...,0,\zeta)$, $\zeta \in \cc$ and $E = \{ r^j(Z) = 0,
j = 1,...,n-1 \}$, where $r^j = z_j + \overline z_j + O(\vert Z
\vert^2)$.

Consider a
$J$-holomorphic disc tangent to  $H_0^J(E)$ at the center.
Performing if necessary an appropriate
 diffeomorphism with   linear part identity at the origin we  can assume
that this disc lies on  $l$.
Thus, we shall assume that $l$ is $J$-holomorphic.

\begin{e-lemme}
In the above variables,
 for every $j$ the Levi form
$L_0^J(r^j)$  coincides on $ H_0^J(E)$ with the Levi form  $L_0^{J_{st}}(r^j)$
 with respect to $J_{st}$.
\end{e-lemme}

\proof This follows from Proposition \ref{PROP1} if in its setting
we take the line
$l$ for  a $J$-holomorphic disc $f$.

  For $\delta > 0$ consider now the {\it non-isotropic}
dilations
$\Lambda_{\delta}:(z,w) \mapsto (\delta^{-1} z,\delta^{-1/2}w)$
and the induced structures $J_{\delta}:=(\Lambda_{\delta})_*(J)$.

\begin{e-lemme}
For any positive real $k$ one has $\parallel J_{\delta} -
  J_{st}\parallel_{\CC^k(K)} \longrightarrow 0$  as
  $\delta \longrightarrow 0$ on  each compact subset $K$ of $\cc^n$.
\end{e-lemme}

\proof Consider the Taylor expansion of $J(Z)$ near the origin: $J(Z)
= J_{st} + L(Z) + R(Z)$ where $L(Z)$ is  the linear part of the expansion
 and $R(Z) =
O(\vert Z \vert^2)$. Clearly, $\Lambda_{\delta} \circ
R(\Lambda_{\delta}^{-1}(Z)) \circ \Lambda_{\delta}^{-1}$ converges to $0$
as $\delta \longrightarrow 0$. Denote by $L_{kj}^{\delta}(Z)$
(respectively, by $L_{kj}(Z)$) an entry
of the real  matrix $\Lambda_{\delta} \circ
L(\Lambda_{\delta}^{-1}(Z)) \circ \Lambda_{\delta}^{-1}$
(respectively, of $L(Z)$). Then
$L_{kj}^{\delta}(z,w) = L_{kj}(\delta z, \delta^{1/2}w)
\longrightarrow 0$ for $k,j =
1,...,2n-2$ and $k,j = 2n-1,2n$,  $L_{kj}^{\delta}(z,w) =
\delta^{1/2}L_{kj}(\delta z, \delta^{1/2}w) \longrightarrow 0$
for $k= 2n-1,2n$, $j = 1,...,2n-2$. For $k = 1,...,2n-2$
and $j= 2n-1,2n$
we have $J_{kj}^{\delta}(z,w) = \delta^{-1/2}L_{kj}(\delta
z,\delta^{1/2}w) \longrightarrow L_{kj}(0,w)$. However, in the
 coordinate system fixed above the line $l$ is $J$-holomorphic,
 that is, $J(l(\zeta)) \circ dl = dl \circ J_{st}$. This shows that
$L(0,w) \equiv 0$. Thus,
$L_{kj}^{\delta}$ approaches  $0$ for all $k,j$. This gives us the
result of the lemma.

We point out that this result fails for $CRdim E >
1$. For this reason we begin our construction  with the case
      $CRdim E = 1$.

We may assume that $E$ is defined by  equations $r^j(z,w) = 0$,
$j=1,...,n-1$ with  $r^j(z, w) = 2Re z_j + 2Re Q^j(z,w) + H^j(z,w) + O(\vert
Z \vert^2)$. Here $Q^j(Z) = \sum q_{ks}^jZ_kZ_s$ and $H^j(Z) =
\sum_{ks}h^j_{ks}Z_k\overline Z_s$ are complex and Hermitian quadratic forms,
respectively. Then the manifold $E_{\delta}:= \Lambda_{\delta}(E)$ is given by the
equations $r^j_{\delta}(Z):= \delta^{-1}r_j((\delta^{1/2})z,\delta
w) = 0$ and $r^j_{\delta}(Z) \longrightarrow r^j_0(Z):=  2Re z_j +
2Re Q^j(0,w) + H^j(0,w)$ (in the $C^k$ norm for any $k$) as
$\delta$ approaches $0$. Since the quadratic map
$$(H^1('0,w),...,H^{n-1}('0,w))$$ can be identified with the Levi
form of $E$ at the origin,
one of the forms $H^j(0,\bullet)$ does not vanish on $\cc$. Replacing
the functions
$r^j$ by their linear combinations if necessary one can assume that
$H^j(0,w) \equiv 0$, $j=1,...,n-2$ and $H_j(0,w) = -\vert w
\vert^2$.

 Consider the limit manifold $E_0 = \{ r^j(Z) = 0, j=1,...,m \}$.
After a biholomorphic (with respect to $J_{st}$) change of the
 variables  $(z,w) \mapsto (z',w') = (z + Q(0,w),w)$
(here $Q = (Q_1,...,Q_{n-1})$) we obtain a   manifold $E_0'$ defined
by the equations $Re z_j = 0$, $j = 1,...,n-2$, $ 2 Re z_{n-1}
 = \vert w \vert^2$ (we drop the primes).

Following Boggess-Pitts \cite{BoPi} we consider now the   family
$f:\zeta \mapsto (z(\zeta),w(\zeta)$ of holomorphic
Bishop discs attached to $E_0'$ and defined by the formulae

\begin{eqnarray*}
& &z_j(\zeta) = iy_j, j = 1,...,n-2,\\
& &z_{n-1}(\zeta) = (1/2)\left ( c\overline c + \frac{t^2}{(1 +
    \lambda)^2}(\lambda^2 +1) \right ) + \frac{t\lambda}{ 1 +
    \lambda}\overline c + iy_{n-1} + \left ( \frac{t\overline c}{1
    + \lambda} + \frac{t^2\lambda}{(1 + \lambda)^2} \right ) \zeta,\\
& &w(\zeta) = c + \frac{t(\lambda + \zeta)}{1 + \lambda}
\end{eqnarray*}

This family depends on parameters $y = (y_1,...,y_{n-1})$ ranging in some
neighborhood of the origin in $\R^{n-1}$, real parameters $t > 0$ and
$\lambda \in [0,1]$, and a complex parameter $c$ ranging
in a neighborhood of
the origin in $\cc$. We shall write $f(t,\lambda,y,c)(\bullet)$ for
 discs in this family. We are interested in the maps
 $f(t,\lambda,y,c)(-\lambda)$.
They have the following properties:
\begin{itemize}
\item[(a)] for any $t> 0$ one has $\lim_{\lambda \longrightarrow 1}
  f(t,\lambda,y,c)(-\lambda) = (iy_1,...,iy_{n-2}, (1/2)c\overline
  c + iy_{n-1},c)$, so that the points
$f(t,1,y,c)(-1)$ fill a neigborhood of the origin in
   $E_0'$ as   $(y,c)$ ranges over a neighborhood of the origin in
  $\R^{n-1} \times \cc$ and the corresponding map is a diffeomorphism.
\item[(b)] for any fixed $t > 0$ the  differential of the map
$(\lambda,y,c) \mapsto f(t,\lambda,y,c)(-\lambda)$ evaluated at
  $(1,0,0)$ has the maximum possible rank $n+2$.
 \end{itemize}

Now fix sufficiently small positive $t = t_0 > 0$.  By Theorem
 \ref{theosolvingBishop}
for small $\delta > 0$  there exist
$J_{\delta}$-holomorphic discs $F_{\delta}(\lambda,y,\omega)(\bullet)$ $C^k$-
smoothly depending on $\delta,y,\lambda,\omega$ such that
$F_{\delta}(\lambda,y,\omega)(b\D) \subset E_{\delta}$ and
$F_0(\lambda,y,\omega) = f(t_0,\lambda,y,\omega)$. It follows by
continuity from (a) and  (b) that the range $\tilde E_\delta$ of the map
$(\lambda,y,\omega)
\mapsto F_{\delta}(\lambda,y,\omega)(-\lambda)$ considered for
$\lambda$ close to $1$  is an $n+2$-manifold with boundary  that is the range
of the map $(y,\omega)
\mapsto F_{\delta}(1,y,\omega)(-1)$ and therefore lies in
$ E_{\delta}$. Since this map is
close to $f(t,1,y,\omega)(-1)$ and so has the maximum possible rank $n+1$, its range is entire $ E_{\delta}$.

{\bf Remark.} Our proof allows one to `control' in a certain measure
 the direction in which the
manifold $\tilde E$ is attached to $E$. Indeed, differentiating the
map $f(t,\lambda,y,\omega)(-\lambda)$ with respect to $\lambda$ at the
point $(t, 1, 0,0)$ we see that the tangent space to $\tilde E_0$ at the
origin is spanned by $T_0(E'_0)$ and the vector $\nu=(0,\dots,1,0)$.
Hence the tangent space to $\tilde E$ at the origin is spanned by $T_0(E)$
and a vector close to $\nu$.

\subsection{The case $CRdim E > 1$}



Let $E$ be a generating submanifold in an almost complex manifold
$(M,J)$. In this section we are particularly interesting in the case
 $CR dim E > 1$, but our considerations are also meaningful for
$CR dim E =1$. As before,  we
assume that $M$ is a neighborhood of the origin in $\cc^n$, $J$ is a
smooth matrix valued function,  $J = J_{st} + O(\vert Z
\vert)$,  and $E = \{ r^j(Z) = 0,
j = 1,...,m \}$, where $r^j = z_j + \overline z_j + O(\vert Z
\vert^2)$, $ Z = (z,w) \in \cc^m \times \cc^{n-m}$.

Let $v \in H_0^J(E)$ be  a vector such that the Levi form of $r^m$ does
not vanish on $v$. Consider a
$J$-holomorphic disc tangent to  $v$ at the center. After a
suitable diffeomorphism with  linear part at the origin
that is  $\cc$-linear this disc coincides with the line
$l = (0,...,0,\zeta)$, $\zeta \in \cc$; pushing forward $J$, we still obtain an
almost complex structure coinciding with the standard one at the
origin. Thus, we may assume that $l$ is $J$-holomorphic in our coordinates.
 Similarly to  the previous section, for every defining function $r^j$ the
value of the Levi form
$L_0^J(r^j)(v)$  coincides with that of the Levi form  $L_0^{J_{st}}(r^j)(v)$
 with respect to $J_{st}$  in the above coordinates.

  For $\delta > 0$ consider the
dilations
$\Lambda_{\delta}:(z,w) \mapsto (\delta^{-1} z,\delta^{-1/2}w)$
and the induced structure $J_{\delta}:=(\Lambda_{\delta})_*(J)$. As we
shall see, if the CR dimension of $E$ is $> 1$, then the structures
$J_\delta$ do not converge to $J_{st}$ in general. Consider the
Taylor expansion of the matrix function $J$:
$$
J(Z) = J_{st} + L(Z) + O(\vert Z \vert^2)
$$
where $L(Z)$ is  the linear part. We observe that $L$ is an
endomorphism of $\R^{2n}$ antilinear with respect to the standard
complex structure.  We regard  $L(Z)$ as a complex $n \times
n$-matrix with entries  $L_{qj}(Z)$ that  are $\R$-linear  ($\cc$-valued)
 functions
of $Z$. The following result can be proved by direct computation.

\begin{e-lemme}
On has  $J_\delta \longrightarrow J_0$ as $\delta \longrightarrow
0$, where $J_0 = J_{st} + L_0(w)$ and the matrix $L_0(w)$ (in the
complex notation)  has entries $L^0_{qj}$  described
as follows: $L_{qj}^0 = 0$, for  $q=1,...,n, j=1,...,m$ and for
$q=n-m+1,...,n$, $j=1,...,n$; for $q=1,...,m$, $j= n-m+1,...,n$ one has
$L_{qj}^0(w) = L_{qj}(0,w)$.
\end{e-lemme}

Moreover, the above condition of the $J$-holomorphy of the line $l$
implies that $L_{q,n-m}$ does not depend on $w_{n-m}$, that is,
$L_{q,n-m}^0 = L_{q,n-m}^0(0,w_1,...,w_{n-m-1})$.

Consider the manifolds $E_\delta: = \Lambda_{\delta}(E)$ defined by the
equations $r^J_\delta(z,w):= \delta^{-1}r^j(\delta z,\delta^{1/2}w) =
0$, $j=1,...,m$. Consider the Taylor expansion $r^j(z,w) = z_j +
\overline z_j + 2Re Q_j(z,w) + H_j(z,w) + O(\vert Z \vert^2)$, where
$Q_j$ is  the  complex quadratic part and $H_j$  the
Hermitian part of the expansion. As $\delta \longrightarrow 0$, we
have $r^j \longrightarrow r^j_0:= z_j + \overline z_j + 2Re Q_j(0,w) +
H_j(0,w)$. We point out that the  biholomorphic  (with respect
to $J_{st}$) change of the variables  $(z,w) \mapsto (z + Q(0,w),w)$
(where $Q = (Q_1,...,Q_m)$) does not change the line $l$, therefore we can
execute it before the dilation. This allows us to assume that $ Q_j(0,w)
\equiv 0$. Thus, the functions $r^j_\delta$ converge to $z_j + \overline
z_j + H_j(0,w)$ as $\delta \longrightarrow 0$. In this sense
we view  the manifold
$E_0 = \{ z_j + \overline z_j + H_j(0,w) = 0, j = 1,...,m \}$ as the
limit of  $E_\delta$ as $\delta \longrightarrow 0$.

Our next aim  is the description of  $J_0$-holomorphic
 Bishop discs (with values in a sufficiently small
neighborhood $U$ of the origin)  with
boundaries attached to $E_0$.  Let $f:\D \longrightarrow U$ be a
 smooth map. To simplify the notations, we will denote by $f_\zeta$
 the partial derivative $\frac{\partial f}{\partial \zeta}$.

 Recall that the $J_0$-holomorphy condition for $f$ can be
 written in the following form:

$$
f_{\overline\zeta} + Q(f)\overline{f_\zeta} = 0,
$$
where $Q(Z)$ is the
complex $n\times n$ matrix  of an operator the composite of which with complex
conjugation is equal to the endomorphism
$- (J_{st} + J_\delta(Z))^{-1}(J_{st} - J_\delta(Z))$
(which is an anti-linear operator with respect to the standrard
structure $J_{st}$). If $f$ has the form $f(\zeta) =
(z(\zeta),w(\zeta))$, then after direct computations of the matrix $Q$
we obtain the equations  of the $J_0$-holomorphy of $f$:

\begin{eqnarray}
\label{***}
(z_j)_{\overline\zeta} = -(i/2)\left (\sum_{q=1}^{n-m} L_{jq}^0(w)
\overline{(w_j)_{\zeta}} \right ), j=1,...,m
\end{eqnarray}
and
$$(w_j)_{\overline\zeta} =0, j=1,...,n-m.$$
This gives one  a direct description of all $J_0$-holomorphic discs. Fix a
function $w \in ({\cal O}(\D))^{n-m} \times (C^k(\overline \D))^{n-m}$
  for  some fixed non-integer $k > 0$. Then  integration of the above
  system (\ref{***}) shows  that

$$z_j(\zeta) = -T_{CG}\left ((i/2)\sum_{q=1}^{n-m} L_{jq}^0(w)
\overline{(w_q)_{\zeta}}\right )(\zeta) + \phi_j(\zeta), j=1,..m$$
where, as before,  $T_{CG}$ is the Cauchy-Green transform and
$\phi_j$ is a holomorphic function of the class
$C^k(\overline\D)$.

For  $w  \in ({\cal O}(\D))^{n-m}
\times (C^k(\overline \D))^{n-m}$ let  $\Psi_j(w)(\bullet)$
be the function
$$
-T_{CG}\left ((i/2)\sum_{q=1}^{n-m} L_{jq}^0(w)
\overline{(w_q)_{\zeta}}\right )
$$
If $w$ is fixed, then the boundary condition
$(z,w)(b\D) \subset E_0$ holds  if and only if

$$
z_j(\zeta) = (\Psi_j(w) - I_{S}(Re \Psi_j(w)))(\zeta) -(1/2)
I_S(H_j(0,w))(\zeta) + iy_j, j=1,...,m,
$$
where $y_j \in \R$ and $I_S$ is  the Schwarz integral in the unit
disc:
\begin{eqnarray}
I_S(h)(\zeta) = \frac{1}{2\pi}\int_0^{2\pi}  h(e^{i\tau})\frac{e^{i\tau} +
  \zeta}{e^{i\tau} - \zeta}d\tau,
\end{eqnarray}
This gives us a complete description of Bishop   discs attached to
$E_0$. In particular, we have the following result.
\begin{e-lemme}
A map $(z,w): \zeta \mapsto (z(\zeta),w(\zeta))$ is a
$J_0$-holomorphic Bishop disc
for $E_0$ if and only if $(z - \Psi_j(w) + I_{S}(Re \Psi_j(w)),w)$ is
a $J_{st}$-holomorphic Bishop disc for $E_0$.
\end{e-lemme}

Similarly to the previous subsection, consider the map $w:\D
\longrightarrow \cc^{n-m}$ of the following form:

\begin{eqnarray*}
& &w_1 = c_1,\\
& &..........\\
& &w_{n-m-1} = c_{n-m-1},\\
& &w_{n-m}(\zeta) = c_{n-m} + \frac{t(\lambda + \zeta)}{1 + \lambda}
\end{eqnarray*}
where the $c_j$ are complex constants, $t> 0$ and $\lambda \in [0,1]$.
Then $(w_k)_\zeta = 0$ for $k = 1,...,n-m-1$. On the other hand, by
our construction $L_{jn-m}^0(0,...,0,\zeta) \equiv 0$. Recall here
that the $L_{jk}^0$ participate in the linear part $L(Z)$ of the Taylor expansion
of $J$ at the origin and are $\R$-linear in $w$. Hence  the
$\R$-linear  function $L_{jn-m}^0$ is independent of $w_{n-m}$ and
$\sum_{k=1}^{n-m} L_{jk}^0(w)\overline{(w_k)_{\zeta}} =
L_{jn-m}^0(c_1,...,c_{n-m-1})t/(1+\lambda)$ is  constant with respect
to the variable $\zeta$.

Let  $a_j(c_1,...,c_{n-m-1}) = -
(i/2)L_{jn-m}^0(c_1,...,c_{n-m-1})t/(1+\lambda)$. Then it
follows from   equations
(\ref{***}) that  $z_j(\zeta) = a_j\overline \zeta +\phi_j(\zeta)$ with
holomorphic $\phi_j$. We set $\phi_j = - \overline a_j  \zeta +
\Phi_j$. Then the inclusion $(z,w)(b\D) \subset E_0$ is equivalent to
the relation $Re \Phi_j (\zeta) = -(1/2)H_j(0,w(\zeta))$ for $\zeta
\in b\D$ meaning that $(\Phi,w)$ is a Bishop disc for $E_0$ with
respect to $J_{st}$.  In view of the condition
$\zeta \overline{\zeta} = 1$, the right-hand
side represents a real polynomial of degree 1 in $\zeta$, so that the $\Phi_j$
are complex polynomials of degree at most 1  and can easily  be
explicitely written; since the  Hermitian quadratic form $H_m(0,w)$
does not vanish on the line $l = (0,...,0,\zeta)$, it contains the
term  with negative coefficient, and
without loss of generality we can assume that this coefficient is $-1$, 
  other terms of the  form $H_m(0,w)$ are independent of $ w_{n-m}$, and
the other forms $H_j(0,w)$ contain no  term $\vert w_{n-m}\vert^2$.
Then we obtain  a formula  for $z_{n-m}$
similar to the one  for $z_{n-1}$ in the previous section and
 explicit
expressions for the $z$-component of a   $J_0$-holomorphic Bishop
disc $(z,w)(c_1,...,c_{n-m},t,\lambda,y_1,...,y_m)(\bullet)$ determined
by the parameters $c_j$, $t$, $\lambda$, $y_k$:

\begin{eqnarray*}
& &z_j(\zeta) = iy_j + a_j(c)\overline{\zeta} - \overline{a}_j(c)\zeta
-\frac12H_j\left (
0,c_1,...,c_{n-m-1},c_{n-m}+\frac{t\lambda}{1+\lambda}\right )- 
\frac{t}{1+\lambda}\zeta l_j(c),\\
& &j = 1,...,m-1,\\
& &z_{n-m}(\zeta) = (1/2)\left ( c_{n-m}\overline c_{n-m} + \frac{t^2}{(1 +
    \lambda)^2}(\lambda^2 +1) \right ) + \frac{t\lambda}{ 1 +
    \lambda}\overline c_{n-m} + iy_{n-1} +\\
& & \left ( \frac{t\overline c_{n-m}}{1
    + \lambda} + \frac{t^2\lambda}{(1 + \lambda)^2} \right ) \zeta
+ a_m(c)\overline{\zeta} - \overline{a}_m(c) \zeta-
\frac12H_{n-m}(0,c_1,...,c_{n-m-1},0),
\end{eqnarray*}
where the $a_j(c)$ are defined above and the $l_j$
are homogeneous linear forms of $c_1,...,c_{n-m-1}$.
 As pointed out already,
 for $c_j =0$, $j = 1,...,n-m+1$ these  are just
$J_{st}$-holomorphic Bishop discs.

Finally, it   is easy to see (by computing the rank of the corresponding
map; cf.   the previous subsection) that the constructed family of
$J_0$-holomorphic Bishop   discs sweeps out a manifold with  boundary
$E_0$. So we may use   the implicit function
theorem to construct a  perturbed family of $J_\delta$-holomorphic
Bishop discs sweeping out a manifold with  boundary $E_\delta$.

\end{document}